\documentclass[14pt, twoside, a4paper]{article}
\let\getprepared\relax
\ifx\TestIngCommand\undefined\relax
\else\input ../../../proc/ppreamb-book.tex 
  \input init.tex \fi
\let\TestIngCommand\undefined

\usepackage{amsmath,amscd,amssymb}
\usepackage{graphics} 
\newtheorem{theo}{Theorem}
\newtheorem{lem}{Lemma} 
\newtheorem{cor}{Corollary} 
 
\newtheorem{defi}{Definition}

\newskip\ttglue\ttglue=.5em plus.25em minus.15em 
\def\firstname#1{\def\FIRSTNAME{#1}\ignorespaces}
\def\lastname#1{\def\LASTNAME{#1}\ignorespaces}
\def\middleinitial#1{\def\MIDDLEINI{#1}\ignorespaces}
\def\department#1{\def\DEPARTMENT{#1}\ignorespaces}
\def\institute#1{\def\INSTITUTE{#1}\ignorespaces}
\def\address#1{\def\ADDRESS{#1}\ignorespaces}
\def\country#1{\def\COUNTRY{#1}\ignorespaces}
\def\otheraffiliation#1{\def\OTHERAFFILIATION{#1}\ignorespaces}
\def\email#1{\def\EMAIL{#1}\ignorespaces}
\newcount\autcount\autcount=0 
\newcount\affcount\affcount=0 
\newcount\numcount\newcount\nummcount 
\newcount\nummmcount\newcount\nummmmcount 
\def\writename#1#2{\ \kern-1ex\hbox{
  \csname AUthor\the#1\endcsname\
  \edef\TESTSTR{}\expandafter\ifx\csname auTHor\the#1\endcsname\TESTSTR
  \else\csname auTHor\the#1\endcsname.\ \fi 
  \csname authOR\the#1\endcsname$^{\csname AFF\the#1\endcsname}$
  \expandafter\ifx\csname corr\number#1\endcsname\relax
  \else\thanks{Corresponding author.}\ \fi 
  }\ifnum#1<#2, \else\ \kern-1ex\fi}
\def\writeemail#1{
  \nummcount=0\relax\nummmcount=0\relax 
  \loop\ifnum\nummcount<\autcount\advance\nummcount by1\relax
    {\expandafter\ifnum\csname AFF\the\nummcount\endcsname=#1\relax
    \global\advance\nummmcount by1\fi}\repeat 
  \nummcount=0\relax\nummmmcount=0\relax 
  \loop\ifnum\nummcount<\autcount\advance\nummcount by1\relax
    {\expandafter\ifnum\csname AFF\the\nummcount\endcsname=#1\relax
    \global\advance\nummmmcount by1\relax\def\blank{}\expandafter 
    \ifx\csname EMAIL\the\nummcount\endcsname\blank(no e-mail)
    \else\csname EMAIL\the\nummcount\endcsname 
    \fi
    \ifnum\nummmmcount<\nummmcount; \fi\fi}\repeat}
\long\def\BeginAuthorList#1\EndAuthorList{#1\relax
  \author{\vbox{\hsize=390pt\noindent\numcount=0\relax
    \loop\ifnum\numcount<\autcount\advance\numcount by1\relax
      \writename{\numcount}{\autcount}
      \repeat}\\[2mm]
    \vbox{\small\numcount=0\relax                                          
      \loop\ifnum\numcount<\affcount\advance\numcount by1\relax            
        \vbox{{\count0=\numcount\relax                                     
          \loop\expandafter\ifnum\csname AFF\the\count0\endcsname
            <\numcount\relax\advance\count0 by1\relax\repeat               
          $^{\csname AFF\the\count0\endcsname}$}
        \def\BLANK{}\expandafter\ifx\csname DEPT\the\numcount\endcsname    
          \BLANK                                                           
          \else\csname DEPT\the\numcount\endcsname, \fi                    
        \csname INST\the\numcount\endcsname,                               
        \csname ADDR\the\numcount\endcsname,                               
        \csname COUN\the\numcount\endcsname                                
        \edef\TEST{}\expandafter\ifx\csname OTHE\the\numcount\endcsname
          \TEST                                                            
          .\else;\break\csname OTHE\the\numcount\endcsname.\fi}
        \vbox{\writeemail{\numcount}}
        \repeat}\\}}
\expandafter\def\csname x1\endcsname{}
\expandafter\def\csname x2\endcsname{}
\expandafter\def\csname x3\endcsname{}
\expandafter\def\csname x4\endcsname{}
\expandafter\def\csname x5\endcsname{}
\expandafter\def\csname x6\endcsname{}
\expandafter\def\csname x7\endcsname{}
\expandafter\def\csname x8\endcsname{}
\expandafter\def\csname x9\endcsname{}
\def\Author#1#2{\global\advance\autcount by1\relax#2                       
  \expandafter\edef\csname AUthor\the\autcount\endcsname{\FIRSTNAME}
  \expandafter\edef\csname auTHor\the\autcount\endcsname{\MIDDLEINI}
  \expandafter\edef\csname authOR\the\autcount\endcsname{\LASTNAME}
  \expandafter\edef\csname EMAIL\the\autcount\endcsname{\EMAIL}
  \let\tempera\"\def\"{\string\"}\expandafter\ifx\csname x\DEPARTMENT
    \endcsname\relax                                                       
    \global\advance\affcount by1\relax\let\"\tempera                       
    \expandafter\edef\csname DEPT\the\affcount\endcsname{\DEPARTMENT}
    \expandafter\edef\csname INST\the\affcount\endcsname{\INSTITUTE}
    \expandafter\edef\csname ADDR\the\affcount\endcsname{\ADDRESS}
    \expandafter\edef\csname COUN\the\affcount\endcsname{\COUNTRY}
    \expandafter\edef\csname OTHE\the\affcount\endcsname{\OTHERAFFILIATION}
    \expandafter\edef\csname AFF\the\autcount\endcsname{\the\affcount}
  \else\expandafter\edef\csname AFF\the\autcount\endcsname{\DEPARTMENT}
  \fi\let\"\tempera\ignorespaces}
\def\CorrespondingAuthor#1#2{
  \expandafter\xdef\csname corr\number#1\endcsname{cor}
  \Author#1{#2}}
\def\PaperTitle#1{\title{\bf#1}}
\def\Category#1{\ignorespaces}
\def\keywords#1{{\noindent \emph{Keywords:}                                
  \def\BLANK{}\def\TEST{#1}\ifx\BLANK\TEST(n/a).\else#1\fi}}
\getprepared
\begin{document}

\PaperTitle{Tutorial on Rational Rotation $C^*$--Algebras}

\Category{(Pure) Mathematics}

\date{}

\BeginAuthorList
  \Author1{
    \firstname{Wayne}
    \lastname{Lawton}
    \middleinitial{M}
    \department{Department of the Theory of Functions, Institute of Mathematics and Computer Science}
    \institute{Siberian Federal University}
    \otheraffiliation{}
    \address{Krasnoyarsk}
    \country{Russian Federation}
    \email{wlawton50@gmail.com}}
\EndAuthorList
\maketitle
\thispagestyle{empty}
\begin{abstract}
The rotation algebra $\mathcal A_{\theta}$ is the universal $C^*$--algebra generated by unitary operators $U, V$ satisfying the commutation relation $UV = \omega V U$ where $\omega= e^{2\pi i \theta}.$ They are rational if $\theta = p/q$ with $1 \leq p \leq q-1,$ othewise irrational. Operators in these algebras relate to the quantum Hall effect \cite{boca,rammal,simon}, kicked quantum systems \cite{lawton1, wang}, and the spectacular solution of the Ten Martini problem \cite{avila}. 
Brabanter \cite{brabanter} and Yin \cite{yin} classified rational rotation $C^*$--algebras up to $*$-isomorphism. Stacey \cite{stacey} constructed their automorphism groups. They used methods known to experts: cocycles, crossed products, Dixmier-Douady classes, ergodic actions, K--theory, and Morita equivalence. 
This expository paper defines $\mathcal A_{p/q}$ as a $C^*$--algebra generated by two operators on a Hilbert space and uses linear algebra, Fourier series and the Gelfand-Naimark-Segal construction \cite{gelfand} to prove its universality. It then represents it as the algebra of sections of a matrix algebra bundle over a torus to compute its isomorphism class. The remarks section relates these concepts to general operator algebra theory. We write for mathematicians who are not $C^*$--algebra experts.
\end{abstract}
%
\footnote{\thanks{This work is supported by the Krasnoyarsk Mathematical 
Center and financed by the Ministry of Science and Higher Education 
of the Russian Federation in the framework of the establishment and 
development of Regional Centers for Mathematics Research and 
Education (Agreement No. 075-02-2020-1534/1).}}

\section{Uniqueness of Universal Rational Rotation $C^*$--algebras}\label{Sec1}
$\mathbb{N},$ $\mathbb{Z},$ $\mathbb{Q},$ $\mathbb{R},$ $\mathbb{C}$ and
$\mathbb{T} \subset \mathbb{C}$ denote the sets of positive integer, integer, rational, 
real, complex and unit circle numbers. For a Hilbert space $H$ let $\mathcal B(H)$ be the $C^*$--algebra of bounded operators on $H.$ All homomorphisms are assumed to be continuous. We assume famliarity with the material in Section \ref{Sec4}.
\\ \\
Fix $p, q \in \mathbb N$ with $p \leq q-1$ and gcd($p,q$) $= 1,$ define $\sigma := e^{2\pi i/q}$ and $\omega := \sigma^p,$ and let $\mathfrak C_{p/q}$ 
be the set of all $C^*$--algebras generated by a set $\{U, V\} \subset \mathcal B(H)$ satisfying 
$UV = \omega VU.$
Since $\{U, V\} = \{V, U\},$ 
$\mathfrak C_{(q-p)/q} = \mathfrak C_{(q-p)/q}.$
$M_q$ and the circle subalgebra of $L^2(\mathbb T)$ generated by $(Uf)(z) := zf(z)$ and $(Vf)(z) := f(\omega z)$ belong to $\mathfrak C_{(q-p)/q}.$ The circle algebra is isomorpic to the tensor product $C(\mathbb T) \otimes M_q.$
\begin{defi}\label{universal}
	$\mathcal A \in  \mathfrak C_{p/q}$ generated by 
	$\{U,V\} \subset \mathcal B(H)$ satisfying 
	$UV = \omega VU$ is called universal if for every 
	$\mathcal  A_1 \in  \mathfrak C_{p/q}$ generated 
	by $\{U_1,V_1\} \subset \mathcal B(H_1)$ satisfying 
	$U_1V_1 = \omega V_1U_1,$ 
	there exists a $*$-homomorphism 
	$\Psi : \mathcal A \mapsto \mathcal A_1$ 
	satisfying $\Psi(U) = U_1$ and $\Psi(V) = V_1.$
\end{defi}
\begin{lem}\label{unique}
		If $\mathcal A, \mathcal A_1 \in  \mathfrak C_{p/q}$ are both universal, then they are isomorphic.
\end{lem}
{\bf Proof}
Let $U,V,U_1,V_1$ be as in Definition \ref{universal}. There exists $*$--homomorphisms $\Psi : \mathcal A \mapsto \mathcal A_1$
and $\Psi_1 : \mathcal A_1 \mapsto \mathcal A$ with
$\Psi_1 \circ \Psi(U) = U,$ $\Psi_1 \circ \Psi(V) = V,$
$\Psi \circ \Psi_1(U_1) = U_1,$ $\Psi \circ \Psi_1(V_1) = V_1.$
Since $\{U,V\}$ generates $\mathcal A,$ $\Psi_1 \circ \Psi$ is the identity map on $\mathcal A.$ Similarly, $\Psi \circ \Psi_1$ is the identity map on $\mathcal A_1.$ Therefore $\Psi$ is a $*$--isomorphism of $\mathcal A$ onto $\mathcal A_1$ and $A$ is $*$--isomorphic to $A_1.$
\section{Construction of Universal Rational Rotation $C^*$--algebras}\label{Sec2}
Define the Hilbert space $H_q := L^2(\mathbb R^2,\mathbb C^q)$
consisting of Lebesgue measurable $v : \mathbb R^2 \mapsto \mathbb C^q$ satisfying
$
\int_{\mathbb R^2} v^*v \, < \infty,
$
equipped with the scalar product 
$
<v,w> \, := \int_{\mathbb R^2} w^*v.
$
Define $\mathcal P_q$ to be the subset of continuous 
$a : \mathbb R^2 \mapsto \mathcal M_q$  satisfying
\begin{equation}\label{Pq}
	 a(x_1,x_2) = a(x_1+q,x_2) = a(x_1,x_2+q), \ (x_1,x_2) \in \mathbb R^2
\end{equation}
and regarded as a $C^*$--subalgebra of $\mathcal B(H_q)$ acting by
$
	(av)(x) := a(x)v(x), \ \ a \in \mathcal P_q, \, v \in H_q.
$
The operator norm of $a \in \mathcal P_q$ satisfies
\begin{equation}\label{operatornorm}
||a|| = \max_{x \in [0,q]^2} ||a(x)||. 
\end{equation}
Define $U, V \in \mathcal P_q$ by
\begin{equation}\label{UV}
U(x_1,x_2) := e^{2\pi i x_1/q} U_0, \ \ V(x_1,x_2) := e^{2\pi i x_2/q} V_0,
\end{equation}
where $U_0, V_0 \in \mathcal M_q$ are defined by (\ref{U0V0}),
and define $\mathcal A_{p/q}$ to be the $C^*$--subalgebra of $\mathcal P_q$ generated by $\{U,V\}.$ 
Choose $r \in \{1,...,q-1\}$ such that $pr = 1 \mod q.$ Then $r$ is unique, gcd($r,q$) $= 1.$ Define $\sigma := e^{2\pi i/q}$ and 
$\omega := \omega^p.$ Then $\omega^r = \sigma.$ 
\begin{theo}\label{Apqcharacterzation} If $a \in \mathcal A_{p/q}$ then 
\begin{equation}\label{conjV0U0}
	a(x_1+1,x_2) = V_0^{-r}a(x_1,x_2)V_0^r \hbox{ and }  a(x_1,x_2+1) = U_0^ra(x_1,x_2)U_0^{-1}.
\end{equation}
Conversely, if $a \in \mathcal P_q$ satisfies (\ref{conjV0U0}), then $a \in \mathcal A_{p/q}.$
\end{theo}
{\bf Proof} (\ref{UV}) and (\ref{U0V0commutation}) give
 $V^{-r}UV^r = \sigma U$ and $U^rVU^{-r} = \sigma V.$ If $a = U^mV^n,$
then
$$a(x_1+1,x_2) = \sigma^m a(x_1,x_2) = V_0^{-r}a(x_1,x_2)V_0^r; \, 
a(x_1,x_2+1) = \sigma^n a(x_1,x_2) = U_0^{r}a(x_1,x_2)U_0^{-r}.$$
The first assertion follows since span$\{U^mV^n:m,n \in \mathbb Z\}$ is dense in $\mathcal A_{p/q}.$ Conversely, if $a \in \mathcal P_q,$ then (\ref{Pq}),
Lemma \ref{U0V0basis}, and Weierstrass' approximation theorem implies that there exist unique
$c(m,n,j,k) \in \mathbb C$ with
$$
	a(x_1,x_2) \sim \sum_{(m,n) \in \mathbb Z^2} \sum_{j,k = 0}^{q-1}
	c(m,n,j,k)\, e^{2\pi i(mx_1+nx_2)/q} \, U_0^jV_0^k
$$
where $\sim$ denotes Fourier series.
Then (\ref{conjV0U0}) gives
$
	c(m,n,j,k)\sigma^m = c(m,n,j,k) \sigma^j
$
and
$
	c(m,n,j,k)\sigma^n = c(m,n,j,k) \sigma^k.
$ 
Since $\sigma^q = 1,$ $c(m,n,j,k) = 0$ unless $j = m \mod q$and $k = n \mod q.$ 
Define $c(m,n) := c(m,n,m \mod q,n \mod q).$
Then $a \in A_{p/q}$ since
$$
	a \sim \sum_{(m,n) \in \mathbb Z^2} c(m,n) U^mV^n.
$$
Representations $\rho_1, \rho_1 : \mathcal A \mapsto \mathcal B(H)$ 
of a $C^*$--algebra $\mathcal A$ are unitarily equivalent if there exists $U \in \mathcal U(H)$ such that
%
$\rho_2(a) = U\rho_1(a)U^{-1}, \ \ a \in \mathcal A.$
%
%
\begin{theo}\label{irr} If $\mathcal A \in \mathfrak C_{p/q}$ is generated by $\{U,V\}$ with $UV = \omega VU$ and $\rho : \mathcal A \rightarrow \mathcal B(H)$ is an irreducible representation then:
\begin{enumerate}
\item dim $H = q$ so $B(H) = \mathcal M_q,$
\item there exist $z_1, z_2 \in \mathbb T$ such that 
$\rho = \rho_{z_1,z_2}$ where
%
$\rho_{z_1,z_2}(U^jV^k) := z_1^jz_2^kU_0^jV_0^k.$
%
\item $\rho_{z_1^{\prime},z_2^{\prime}}$ is unitarity equivalent to $\rho_{z_1,z_2}$
iff $(z_1^{\prime}/z_1)^q = (z_2^{\prime}/z_2)^q = 1.$
\end{enumerate}
\end{theo}
{\bf Proof} Boca gives a proof in ([1], p. 5, Lemma 1.8, p. 7, Theorem 1.9). We give a proof based on Schur's lemma. Let $\mathcal C \subset \mathcal A$ be the $C^*$--subalgebra generated by $\{U^q, V^q\}.$ Since $\rho$ is irreducible and $\rho(\mathcal C)$ commmutes with
$\rho(\mathcal A),$ 
there exists a $*$--homomorphism $\gamma :  \mathcal C \mapsto \mathbb C$ such that 
$\rho(c) = \gamma(c) I, \ c \in \mathcal C.$
Choose $h \in H\backslash \{0\}$ and define $H_1 := \hbox{span }\{\rho(U^jV^k)h; 0 \leq j, k \leq q-1\}.$ Since $H_1$ is closed, $\rho$--invariant, $H_1 \neq \{0\},$ and $\rho$ is irreducible,
$H = H_1.$ Since dim $H \leq q^2,$ $\rho(V)$ has an eigenvector $b$ with eigenvalue $\lambda \in \mathbb T$ and $||b|| = 1.$ Define $z_2 := \lambda \omega.$ Choose $z_1 \in \mathbb T$ so $z_1^q = \gamma(U^q)$ and define $b_j := z_1^{j} \rho(U^{-j})b, \, 1 \leq j \leq q.$
Then $\rho(V)b_j = z_2 \omega^{j-1} b_j, \, j = 1,...,q,$ and $\rho(U)b_1 = z_1b_q,$ and 
$\rho(U)b_j = z_1 b_{j-1}, \, 2 \leq j \leq q.$
Therefore $\{b_1,...,b_q\}$ is a basis for $H,$ and (\ref{U0V0}) implies that $\rho(U) = z_1U_0,$ and $\rho(V) = z_2V_0$ with respect to this basis. This proves assertions 1 and 2. Assertion 3 follows since
the set of eigenvalues of $\rho(U)$ is $\{z_1\omega^j,\, 0 \leq j \leq q-1\},$ the
set of eigenvalues of $\rho(V)$ is $\{z_2\omega^j,\, 0 \leq j \leq q-1\},$ and the set of eigenvalues determines unitary equivalence. 
 %
 %
 \begin{theo}\label{Apquniversal}
 $\mathcal A_{p/q} \subset \mathcal B(H)$ is the universal $C^*$--algebra in $\mathfrak C_{p/q}.$ 
 \end{theo}
 {\bf Proof} Assume that $\mathcal B \in \mathfrak C_{p/q}.$ Then there exists a Hilbert space
 $H_1$ and $U_1,V_1 \in \mathcal B(H_1)$ with $U_1V_1 = \omega V_1U_1$ and $\mathcal B$ is generated by $\{U_1,V_1\}.$ It suffices to construct a continuous $*$--homomorphism 
$\varphi : \mathcal A_{p/q} \mapsto \mathcal B$
satisfying $\varphi(U) = U_1$ and $\varphi(V) = V_1.$ Define dense $*$--subalgebras
$$
\widetilde {\mathcal A_{p/q}} := \hbox{ span } \{U^jV^k: j, k \in \mathbb Z\} \subset \mathcal A_{p/q}, \ \ 
\widetilde {\mathcal B} := \hbox{ span } \{U_1^jV_1^k: j, k \in \mathbb Z\} \subset \mathcal B,
$$
and a $*$--homomorphism 
$\widetilde \varphi : \widetilde {\mathcal A_{p/q}} \mapsto \widetilde {\mathcal B}$ by
$\widetilde \varphi (U^jV^k) : = U_1^jV_1^k.$ To extend $\widetilde \varphi$ to $*$--homomorphism $\varphi : \mathcal A_{p/q} \mapsto \mathcal B$ it suffices to show that for every Laurent polynomial of two variables $p(u,v)$ the following inequality is satisfied
%
$||p(U_1,V_1)|| \leq ||p(U,V)||$
%
since $p(U_1V_1) = \widetilde {\varphi} (p(U,V)).$ 
Then (\cite{davidson}, Corollary I.9.11), which follows directly from the Gelfand-Naimark-Segal construction, implies that there exists an irreducible representation
$\rho_1 : \mathcal B \mapsto \mathcal M_q$ and $v \in H_1$ with $||v|| = 1$ such that
%
$||p(U_1,V_1)|| = ||\rho_1(p(U_1,V_1))v||.$
%
Theorem \ref{irr} implies that $\rho_1(U_1) = z_1U_0$ and
$\rho_1(V_1) = z_2V_0$ for some $z_1, z_2 \in \mathbb T.$
Let $\rho : \mathcal A_{p/q} \mapsto \mathcal M_q$ be the 
irreducible representation defined by Theorem \ref{irr} so
$\rho(U) = z_1U_0$ and $\rho(V) = z_2V_0.$
Since $\rho_1 \circ \widetilde \varphi = $ the restriction of $\rho$ to
$\widetilde {\mathcal A_{p/q}},$ (\ref{operatornorm}) and (\ref{UV}) imply that
%
$$||p(U_1,V_1)|| = ||\rho_1(p(U_1,V_1))v|| \leq ||\rho(p(U,V))|| \leq ||p(U,V)||$$
%
which concludes the proof.
\section{Bundle Topology and Isomorphism Classes}\label{Sec3}
Define $\mathbb E_1$ to be the Cartesian product $[0,1]^2 \times \mathcal M_q$
with the identification
%
$$(1,x_2,M) = (0,x_2,V_0^{-r}MV_0^r), \ \ x_2 \in [0,1], M \in \mathcal M_q$$
%
and
%
$$(x_1,1,M) = (x_1,0,U_0^rMU_0^{-r}), \ \ 
x_1 \in [0,1], M \in \mathcal M_q$$
%
and define the algebra bundle $\pi_1 : \mathbb E_1 \mapsto \mathbb T^2$ by
%
$$\pi_1(x_1,x_2,M) = (e^{2\pi i x_1},e^{2 \pi i x_2}), \ \ (x_1,x_2,M) \in \mathbb E_1.$$
%
A map $s : \mathbb T^2 \mapsto \mathbb E_1$ is called a
section if it is continuous and $\pi_1 \circ s = I$ where $I$ denotes the identity map on $\mathbb T^2.$ Since for every
$p \in \mathbb T^2,$ the fiber $\pi_1^{-1}(p) = \mathcal M_q,$
the set of sections under pointwise operations is a $C^*$--algebra. The theorems above show that this algebra is isomorphic to $\mathcal A_{p/q}.$ Furthermore, since points in $\mathbb T^2$ correspond to unitary equivalence classes of irreducible representations, isomorphism of algebras induces homeomorphisms of $\mathbb T^2.$ In order to compute isomorphism classes of universal rational rotation $C^*$--algebras it is convenient to use a slightly different bundle representation of $\mathcal A_{p/q}.$ Define $W \in \mathcal P_q$
%
$$W(x_1,x_2) := 
				\left[
         \begin{array}{cccc}
          1 & 0 & 0 & 0 \\
          0 & e^{2\pi i x_1/q} & 0 & 0 \\     
          0 & 0 & \ddots & 0 \\
          0 & 0 & 0 &e^{2\pi i (q-1)x_1/q}
          \end{array}
          \right]$$
%
and $\mathcal A_{p/q}^{\, \prime} := W \mathcal A_{p/q}W^{-1},$ which is $*$--isomorphic to $\mathcal A_{p/q}.$ Then 
$\mathcal A_{p/q}^{\, \prime}$ is represented as the algebra of sections of the
algebra bundle $\pi_2 : \mathbb E_2 \mapsto \mathbb T^2$ where 
$\mathbb E_2$ is the Cartesian product $\mathbb T \times [0,1] \times \mathcal M_q$
with the identification
%
$$(z_1,1,M) = (z_1,0,G^rMG^{-r}), \ \ z_1 \in \mathbb T, M \in \mathcal M_q$$
%
and 
%
$$G(z_1) := 
				\left[
         \begin{array}{ccccc}
          1 & 0 & 0 & 0 & 0 \\
          0 & 1 & 0 & 0 & 0 \\     
          0 & 0 & \ddots & 0 & 0 \\
          0 & 0 & 0 & 1 & 0 \\
          0 & 0 & 0 & 0 & z_1
          \end{array}
          \right]\, U_0.$$
%
$G^r$ is the clutching function of the bundle.
Let $G_{in} : \mathbb T \mapsto Aut_q^*$ be the map defined by conjugation by $G.$ 
Using the arguments for vector bundles in \cite{hatcher}, it can be shown that the
isomorphism classe of $\mathcal A_{p/q}$ is determined the homotopy class of 
$G_{in}^r : \mathbb T \mapsto Aut_q^*.$ Since $\pi_1(G_{in}) = -1,$ 
$\pi_1(G_{in}^r) = -r$ which gives:
%
%
\begin{theo}\label{isoclass}
$\mathcal A_{p/q}$ is isomorphic to $\mathcal A_{p^{\, \prime}/q^{\, \prime}}$ iff
$q^{\, \prime} = q$ and either $p^{\, \prime} = p$ or $p^{\, \prime} = q-p.$
\end{theo}
\section{Requisite Results}\label{Sec4}
\subsection{Hilbert Spaces and Adjoints}
$H$ is a Hilbert space with
inner product 
$< \cdot, \cdot> : H \times H \mapsto \mathbb C,$
norm
$||v|| := \sqrt {<v,v>},$ and metric $d : H \times H \rightarrow [0,\infty)$ defined by
$d(v,w) : = ||v-w||.$ 
$\mathcal B(H)$ is the Banach algebra of bounded operators on $H$  (continuous linear maps from $H$ to $H$) with operator norm
$$
	||a|| := \sup \{\, ||av||\, :\, v \in H,\, ||v|| = 1 \, \}.
$$
The dual space $H^*$ is the set of continuous linear functions
$L : H \mapsto \mathbb C.$ For $w \in H$ define $L_w \in H^*$ by
$L_w v := \, <v,w>, \ \ v \in H.$ 
\begin{lem}\label{RR}
If $L \in H^{*}$ then there exists a unique $w \in H$ such that $L = L_w.$ 
\end{lem}
{\bf Proof} Rudin gives a direct proof (\cite{rudin}, Theorem 4.12). If $\mathfrak B$ is an orthonormal basis for $H$ and
$
	w := \sum_{b \in \mathfrak B} \overline {Lb}\,  b,
$
then since for every $v \in H,$ 
$v = \sum_{b \in \mathfrak B} <v,b> \,  b,$ it follows that
$$
	L v = 
	\sum_{b \in \mathfrak B} <v,b>\, Lb = 
	\left<v,\, \sum_{b \in \mathfrak B}\overline {Lb} \, \,b\right> \, =\,
	<v,w> \, = L_w \,v.
$$
Lemma \ref{RR} ensures the existence of adjoints. 
For $a \in \mathcal B(H)$ define its adjoint $a^* \in \mathcal B(H)$ by 
$
	L_{a^*w} := L_w \circ a, \ \ w \in H
$
where $\circ$ denotes composition of functions. Therefore
$$
	<a\, \, v,w> = \, <v,a^*w>, \ \ v, w \in H.
$$
Clearly $a^{**} = a,$ $(ab)^* = b^*a^*,$ and the Cauchy--Schwarz inequality gives
$$||a^*|| =  \sup \{\, |<a^*v, w>|:v,w \in H,\, ||v|| = ||w|| = 1\} = $$
$$\sup \{\, |<v, aw>|:v,w \in H,\, ||v|| = ||w|| = 1\} = ||a||$$
and
\begin{equation}\label{Cstar}
	||a^*a|| = \sup \{\, |<a^*av, w>|:v,w \in H,\, ||v|| = ||w|| = 1\} =
\end{equation}
$$\sup \{\, |<av, aw>| :v, w \in H,\, ||v|| = ||w|| = 1 \} = ||a||^2.$$
(\ref{Cstar}) is called the $C^*$--identity. It makes $\mathcal B(H)$ equipped with the adjoint a $C^*$--algebra. The identity operator $I \in \mathcal B(H)$ is defined by $Iv := v$ for all $v \in H.$
$$
	\mathcal U(H) := \{U \in \mathcal B(H):UU^* = U^*U = I\},
$$
the set of unitary operators, is a group under multiplication.
A subalgebra $\mathcal A \subset \mathcal B(H)$ is a 
$C^*$--algebra if it is closed in the metric space topology on $\mathcal B(H)$ and $a^* \in \mathcal A$ whenever $a \in \mathcal A.$ The intersection of any nonempty collection of $C^*$--subalegras of $\mathbb B(H)$ is a $C^*$--algebra. If $S \subset \mathcal B(H)$ the intersection of all $C^*$--subalgebras of $\mathcal B(H)$ that contain $S$ is the $C^*$--algebra generated by $S.$ 
\subsection{Matrix Algebras}
For $m, n \in \mathbb N,$ $\mathbb C^{m \times n}$ 
denotes the set of $m$ by $n$ matrices with complex entries 
and $\mathbb C^n := \mathbb C^{n \times 1}.$ 
The adjoint of $a \in \mathbb C^{m \times n}$ is 
the matrix $a^* \in \mathbb C^{n \times m}$ defined by 
$
a^*_{j,k} := \overline {a_{k,j}}.
$
$\mathbb C^n$ is a Hilbert space with scalar product 
$
	<v,w> \, := \, w^*v, \ \ v, w \in \mathbb C^n.
$
Clearly 
$$
	\mathcal B(\mathbb C^n) = \mathcal M_n
$$
where for $a \in \mathcal M_n$ the adjoint of $a$ as an operator corresponds to the adjoint of $a$ as a matrix. $I_n$ denotes the $n$ by $n$ identity matrix whose diagonal entries equal $1$ and other entries equal $0.$ The operator norm of $a \in \mathcal M_n$ is
$
	||a|| = \sqrt {\hbox{spectral radius } a^*a}
$
where the spectral radius is the largest moduli of the eigenvalues of a matrix.
Thus $\mathcal M_n$ is a $C^*$--algebra. It is also a Hilbert space a Hilbert space of dimension $n^2$ with inner product
\begin{equation}\label{HS}
	<a,b> := \hbox{Trace } b^*a
\end{equation}
and orthonormal basis $e_{j,k} : = $matrix with $1$ in row $j$ and column $j$ with all other enties $= 0.$ 
Fix $p, q \in \mathbb N$ with $p \leq q-1$ and gcd($p,q$) $= 1.$
Define $U_0, V_0 \in \mathcal M_q$ by
\begin{equation}\label{U0V0}
U_0 := \left[
        \begin{array}{cccc}
         0 & 1 &  0 & 0 \\
         \vdots & 0 & \ddots & \vdots \\
         0 & 0 & 0 & 1 \\
         1 & 0 & 0 & 0
         \end{array}
        \right],  \ \ \ \ 
V_0 := 
				\left[
         \begin{array}{cccc}
          1 & 0 & 0 & 0 \\
          0 & \omega & 0 & 0 \\     
          0 & 0 & \ddots & 0 \\
          0 & 0 & 0 & \omega^{q-1}
          \end{array}
          \right],
\end{equation}
\begin{lem}\label{U0V0basis}
$\{(1/\sqrt q)U_0^jV_0^k:0 \leq j, k \leq q-1\}$ is an orthonormal basis for $\mathcal M_q$
with the scalar product defined by (\ref{HS}). Furthermore,
\begin{equation}\label{U0V0commutation}
	U_0V_o = \omega V_0U_0
\end{equation}
\end{lem}
{\bf Proof} (\ref{U0V0commutation}) is obvious. The first assertion  follows since
$$<U_0^jV_0^k,U_0^mV_0^n> \, 
\hbox{Trace }V_0^{-n}U_0^{-m} U_0^jV_0^k = 
\hbox{Trace }U_0^{j-m}V_0^{k-n} = 
\begin{cases}
	q \hbox{ if } j=m \hbox{ and } k=n, \\
	0 \hbox{ otherwise.}
\end{cases}
$$
Define the groups of unitary matrices $\mathcal U_n := \mathcal U(\mathbb C^n)$ and special unitary matrices
$\mathcal S_n := \{a \in \mathcal U_n: \det a = 1\}.$
Clearly $U_0$ and $V_0$ are unitary. Since $\det U_0 = \det V_0 = (-1)^{q-1},$ they are special unitary iff $q$ is odd.
A map $\psi : \mathcal M_n \mapsto \mathcal M_n$ is a homomorphism if
it is linear and satisfies $\psi(ab) = \psi(a)\psi(b)$ 
for all $a, b \in \mathcal M_n$ and an automorphism if is also a bijective. An automorphism $\psi$ is a $*$--automorpism if $\psi(b^*) = \psi(b)^*$ for all $b \in \mathcal M_n.$ 
$Aut_n, Aut_n^*$ 
denote the group of all automorpisms, $*$-automorphisms of $\mathcal M_n.$ 
$\psi \in Aut_n$ is called inner if there 
exists an invertible $a \in \mathcal M_n$ such that 
$\psi(b) = aba^{-1}$ for every $b \in \mathcal M_n.$
%
%
\begin{theo}\label{SN}(Skolem--Noether)
Every $\psi \in Aut_n$ is inner.
\end{theo}
{\bf Proof} The algebra $\mathcal M_n$ is simple, meaning it has no two-sided ideals 
othe that itself (\cite{shilov}, 11.41), so the result follows from the classic Skolem-Noether theorem. An elementary constructive proof is given in \cite{szigeti}. 
%
%
\begin{theo}\label{staraut1}
If $\psi \in Aut_n^*$ then there exists $a \in \mathcal U_n$ such that
$\psi(b) = aba^*$ for every $b \in \mathcal M_n.$
\end{theo}
{\bf Proof} Every $\psi \in Aut_n^*$ induces an irreducible representation $\psi : \mathcal M_n \rightarrow \mathcal B(\mathbb C^n)$ so Theorem \ref{irr} implies that there exists a basis $\{b_1,...,b_n\}$ with respect to which $\psi(U_0)$ has the matrix representation $z_1U_0$ and $\psi(V_0)$ has the representation $z_2V_0.$ Since $U_0^n = V_0^n = I,$
$z_1^n = z_2^n = 1$ so without loss of generality this basis can be chosen to make $z_1=z_2=1$ and then $\psi(a) = aba^{-1}{-1}$ where $ae_j = b_j$ and $\{e_1,...,e_n\}$ is
the standard basis for $\mathbb C^n.$
This theorem can also be derived as a corollary of of Theorem \ref{SN}.
Clearly $\psi(I_n) = I_n.$ Theorem \ref{SN} implies that there exists an invertible $a \in \mathcal M_n$ such that $\psi(b) = aba^{-1}$ for all $b  \in \mathcal M_n.$ Since $\psi$ is a $*$--homomorphism 
$ab^*a^{-1} = (aba^{-1})^* = (a^{-1})^* b^* a^*$
hence $a^*a b^* = b^* a^* a$ for every $b \in \mathcal M_n$ which implies that $a^*a = c\, I_n$ for some $c > 0.$ Replacing $a$ by $a/\sqrt c$ gives the conclusion.
\begin{cor}\label{staraut2}
Let $\mathbb T_n \subset \mathbb T$ be the subgroup of $n$--th roots of 
unity. $\mathbb T_n \, I_n \subset \mathcal S_n$ is isomorphic to 
$\mathbb Z/n\mathbb Z.$ 
$Aut_n^*$ is isomorphic to the quotient group 
$\mathcal S_n/\mathbb T_n \, I_n.$ The fundamental group 
$\pi_1(Aut_n^*)$ is isomorphic to 
$\mathbb Z/n\mathbb Z.$
\end{cor}
{\bf Proof} Assertion one is obvious. Define $\zeta : \mathcal U_n \mapsto Aut_n^*$ by
$\zeta(a)(b) := aba^*.
$
$\zeta$ is a $*$--homomorphism, kernel $\zeta = \mathbb T\, I_n,$ and Corollary \ref{staraut2} implies that $\zeta$ is onto. The first homomorphism theorem of group theory (\cite{waerden}, 7.2) implies that $Aut_n^*$ is isomoprhic to $\mathcal U_n/\mathbb T\, I_n.$ Since $\mathcal n = (\mathbb T I_n)((\mathbb T_n I_n)$ and $\mathbb T_n\, I_n = \mathcal S_n \cap (\mathbb T I_n),$ the second isomorphism theorem of group theory (\cite{waerden}, 7.3) implies that
$Aut_n^*$ is isomporphic to $\mathcal S_n/\mathbb T_n  I_n.$ 
$\mathcal S_n$ is simply connected (\cite{hall}, Proposition 13.11) hence since 
$\mathbb T_n  I_n$ is discrete $\mathcal S_n$ is the univeral cover of $\mathcal S_n \cap (\mathbb T I_n)$ hence the discussion in (\cite{hatcher}, 1.3) implies the last assertion.
\subsection{Spectral Decomposition Theorem for Unitary Operators}
$E \in \mathcal B(H)$ is called a projection if
$E^* =E$ and $E^2 = E.$ Then $P: H \mapsto PH$ is orthogonal projection. A collection of projections $\{E_\varphi: \varphi \in [0,2\pi]$ is called a spectral family if $E_{\varphi_1}E_{\varphi_2} = E_{\varphi_2}E_{\varphi_1} = E_{\varphi_1}$ whenver $\varphi_1 \leq E_{\varphi2}.$ 
\\ \\
Let $A, P, N \in \mathcal B(H).$ $A$ is self--adjoint if $A^* =A.$
$P \in \mathcal B(H)$ is positive if $<Pv,v> \,  \geq \, 0$ for all $v \in H.$
$N \in \mathcal B(H)$ is called normal if $AA^* = A^*A.$
Clearly self-adjoint and unitary operators (or transformations) are normal. Furthermore  eigenvalues of self--adjoint operators are real and eigenvalues of unitary operators have modulus $1.$ If dim $H < \infty,$ then $H$ admits an orthonormal basis of eigenvectors (\cite{shilov}, Theorem 9.33). Therefore every unitary matrix in $\mathcal M_n$ can be diagonalized and its diagonal entries have modulus $1.$ The following result, copied verbatim from the classic textbook by F. Riesz and B. Sz.-Nagy (\cite{riesz}, p. 281),  extends this diagonalization to unitary operators on arbitrary Hilbert spaces.
%
%
\begin{theo}\label{Uspectral}
Every unitary transformation $U$ has a spectral decomposition
$$U = \int_{-0}^{2\pi} e^{i\varphi}dE_{\varphi},$$
where 
$\{E_\varphi \}$ is a spectral family over the segmen 
$0 \leq \varphi \leq 2\pi.$ We can require that 
$E_\varphi$ be continuous at the point 
$\varphi = 0,$ that is, $E_0 = 0;$ 
$\{E_\varphi \}$ 
will then be determined uniquely by $U.$ Moreover, 
$E_\varphi$ is the limit of a sequence of polynomials in $U$ and $U^{-1}.$
\end{theo}
{\bf Proof} The authors of \cite{riesz} reference 1929 papers by von Neumann \cite{neumann} and Wintner \cite{wintner}, 1935 papers by Friedricks and Wecken, and a 1932 book by Stone. They observe that the theorem can be deduced from the one on symmetric transformation (\cite{riesz}, p. 280) (since $U = A + iB$ where
$A := (U+U^*)/2$ and $B := -i(U-U^*)/2$ are symmetric) or from the theorem on trigonometric moments (\cite{riesz}, Section 53), but they give a direct three page proof. We sketch their proof. For every trigomometric polynomial
$p(e^{i\varphi}) = \sum_{-n}^n c_k e^{ik\varphi}$ we associate the transformation $p(U) := \sum_{-n}^n c_k U^k.$ This gives a $*$--homomorphism of the algebra of trigonometric polynomials (where $*$ means complex conjugation) into the subalgebra of $\mathcal B(H)$ generated by $U$ and $U^* = U^{-1}.$ Clearly if $p(e^{i\varphi})$ is real--valued then $p(U)$ is self-adjoint. If $p(e^{i\varphi}) \geq 0$ the Riesz--Fejer factorization lemma (\cite{riesz}, Section 53) implies that there exists a trigonometric polynomial $q(e^{i\varphi})$ with
$
	p(e^{i\varphi}) = q(e^{i\varphi})\overline {q(e^{i\varphi})}
$
hence
$
	p(U) = q(U) q(U)^*.
$
Therefore 
$
	<p(U)v,v> = <q(U)v,q(U)> \, \geq 0, \ \ v \in H,
$
hence $p(U)$ is a positive operator. For $0 \leq \psi \leq 2\pi$ let $e_\psi$ be the characteristic function of $(0,\psi]$ extended to a $2\pi$ periodic function on $\mathbb R.$ Let $p_n$ be a monotonically sequence of positive trigonometric functions with 
$
\lim_{n \rightarrow \infty} p_n(U)v = E_\psi v, \ \ v \in H
$
($p_n(U)$ converges to $E_\psi \in \mathbb B(H)$ in the strong operator topology). $E_\psi$ is a projection since
$E_\psi^* = E_\psi$ and $E_\psi^2 = E_\psi,$ so , and the set $\{E_\varphi:\varphi \in [0,2\pi]$ is a spectral family.
Since the functions $e_\psi$ are upper semi--continuous 
$
	\lim_{\chi \rightarrow \psi, \chi > \psi} E_\chi = E_\psi.
$
Given $\epsilon > 0$ choose $0 < \psi_0 < \psi_1 < \cdots \psi_n = 2\pi$
with max $(\psi_{k+1}-\psi_k) \leq \epsilon$ and choose $\varphi_k \in [\psi_{k-1},\psi_k], k = 1,...,n.$ Then for $\varphi \in (\psi_{k-1},\psi_k]$
$$
	\left|e^{i\varphi} - \sum_{k=1}^n e^{i\varphi_k}\, [e_{\psi_k} - e_{\psi_{k-1}}] \right| = |e^{i\varphi}-e^{i\varphi_k})| \leq |\varphi-\varphi_k| \leq \epsilon
$$
with a similar inequality for $\varphi = 0.$ Since this inequality holds for all $\varphi \in [0,2\pi]$ 
$$
	\left|\left|\,
		U - \sum_{k=1}^n e^{i\varphi_k}\, (E_{\psi_k} - E_{\psi_{k-1}})\,
	\right|\right| \leq \epsilon
$$
A subspace $H_1 \subset H$ is called proper if $H_1 \neq \{0\}$ and $H_1 \neq H.$
The following is an immediate consequence of Theorem \ref{Uspectral}
\begin{cor}\label{corUspectral}
If $U \in \mathcal U(H)$ then either $U = \gamma I$ for some $\gamma \in \mathbb T$ or there exists a projection operator $E : H \mapsto H$ satisfying
\begin{enumerate}
\item $E$ is the limit in the strong operator topology on $\mathcal B(H)$ of polynomials $p(U,U^{-1}).$
\item $EH$ is a proper closed $U$--invariant subspace of $H.$
\end{enumerate}
\end{cor}
\subsection{Irreducible Representations and Schur's Lemma}
A representation of a $C^*$--algebra $\mathcal A$ on a Hilbert space 
$H$ is a $^*$--homomorphism $\rho : \mathcal A \mapsto \mathcal B(H).$ 
A subspace $H_1 \subset H$ is called $\rho$--invariant if 
$\rho(a)H_1 \subset H_1$ for every $a \in \mathcal A.$ $\rho$ 
is irreducible iff it $H$ has no closed proper $\rho$--invariant subspaces. 
The following result extends Shur's lemma for finite dimensional 
representations (\cite{shilov}, 11.33) for unitary operators.
\begin{theo}\label{Schur} (Schur's Lemma)
If $\rho : \mathcal A \mapsto \mathcal B(H)$ is an irreducible representation 
and $U \in \mathcal U(H)$ commutes with $\rho(a)$ for every 
$a \in \mathcal A,$ then there exists $\gamma \in \mathbb T$ with 
$U = \gamma I.$
\end{theo}
{\bf Proof} If the conclusion does not hold then Corollary \ref{corUspectral} 
implies that there exists a projection $E$ satisfying conditions (1) and (2). 
Condition (1) implies that $U\rho(a) = \rho(a)U$ for every 
$a \in \mathcal A.$ The $EH \subset H$ is closed and 
$\rho$--invariant since for every $a \in \mathcal A,$ 
$\rho(a)EH = E\rho(a)H.$ Condition (2) asserts that $EH$ 
is a proper subspace thus contradicting the hypothesis that 
$\rho$ is irreducible, and concluding the proof.
\begin{cor}\label{corSchur}
If $\mathcal A$ is an commutative $C^*$--algebra generated by a set of 
unitary operators and
$\rho : \mathcal A \mapsto \mathcal B(H)$ is irreducible then dim $H = 1$ and
there exists a $*$--homomorphism $\Gamma : \mathcal A \rightarrow \mathcal C$ with $\rho(a) = \Gamma(a)I$ for every $a \in \mathcal A.$
\end{cor}
{\bf Proof} Follows from from Theorem \ref{Schur} since if 
$u \in \mathcal A$ is unitary the $\rho(u) \in \mathcal U(H)$ 
and $\rho(u)$ commutes with 
$\rho(a)$ for every $a \in \mathcal A.$
\section{Remarks} \label{Sec5}
We relate concepts introduced to explain rational rotation algebras to general $C^*$-algebra theory, especially two breakthrough results obtained by teams of computer scientists. 
\\ \\
{\bf Remark 1} Dauns \cite{dauns} initiated a program to represent $C^*$-algebras by continuous sections over bundles over their primitive ideal spaces (kernels of irreducible representations equipped with the hull-kernel topology). The primitive ideal space of rational rotation algebras is homemorpic to the torus $\mathbb T^2.$
\\ \\ 
{\bf Remark 2} Bratteli, Elliot, Evans, and Kishimoto\cite{bratteli} represent fixed point 
$C^*$--subalgebras of $\mathcal A_{p/q}$ by algebras of sections of $M_q$--algebra bundles over the sphere $S^2,$ which is the space of orbits of $\mathbb T^2$ under the map $g \mapsto g^{-1}.$ 
\\ \\
{\bf Remark 3} Elliot and Evans \cite{elliot} derived the structure of irrational rotation algebras. They proved that if $p/q < \theta < p^{\, \prime}/q^{\, \prime},$ then $\mathcal A_{\theta}$ can be approximated by a $C^*$--subalgebra isomorphic to $C(\mathbb T) \otimes M_q \oplus C(\mathbb T) \otimes M_{q^{\, \prime}}.$ This approximation, combined with the continued fraction expansion of $\theta,$ represents $\mathcal A_{\theta}$ as an inductive limit of these subalgebras. 
\\ \\
{\bf Remark 4} Williams \cite{williams} gives an extensive explanation of crossed product $C^*$--algebras, which include rotation algebras.  
\\ \\
{\bf Remark 5} Kadison and Singer \cite{kadison} formulated a problem about extending pure states. Such an extension is used in the Gelfand-Naimark-Segal construction which we used to prove Theorem 3. This problem was shown to be equivalent to numerous problems in functional analysis and signal processing \cite{casazza}, dynamical systems \cite{lawton2,paulsen}, and other fields \cite{bownik}. Weaver \cite{weaver} gave a discrepancy--theoretic formulation that was proved in a seminal paper by three computer scientists: Marcus, Spielman, and Shrivastava \cite{marcus}.
\\ \\
{\bf Remark 6} Courtney \cite{courtney1, courtney2} proved that the class of 
residually finite dimensional $C^*$--algebras, those whose structure can be recovered from their finite dimensional representations, coincides with the class of algebras containing a dense set of elements that attain their norm under a finite dimensional representation, this set is the full algebra iff every irreducible representation is finite dimensional (as for rational rotation algebras), and related these concepts to Conne's embedding conjecture \cite{connes}.  Her publications \cite{courtney3,courtney4,courtney5} cite many references that discuss equivalent formulations of this conjecture.
\\ \\
{\bf Remark 7} In January 2020 five computer scientists: Ji, Natarajan, Vidick, Wright and Yuen submitted a proof that the Conne's embedding conjecture is false. As of November 2021 their paper is still under peer review. However, the editors of the ACM decided, based on the enormous interest that their paper attracted, to publish it \cite{ji}. 
\\ \\
{\bf Acknowledgments} We thank Paolo Bertozzini, Kristin Courtney, S$\phi$ren Eilers, Vern Paulsen, Beth Ruskai, David Yost, and Saeid Zahmatkesh for sharing their insights.


\begin{thebibliography}{10}

\bibitem{avila} A. Avila and S. Jitomirskaya, {\it Solution of the ten martini problem.}
Annals of Mathematics 170, no. 1 (2000) 303--341.

\bibitem{boca} F. Boca, {\it Rotation $C^*$--Algebras and Almost Mathieu Operators.} The Theta Foundation, Bucharest, 2001.

\bibitem{bownik} M. Bownik, {\it The Kadison--Singer problem.} arXiv:1702.04578

\bibitem{brabanter} M. De Brabanter, {\it The structure of rational rotation C¤-algebras.} Arch. Math, 43 (1984) 79--83.

\bibitem{bratteli} O. Bratteli, G. A. Elliot, D. E. Evans, and A. Kishimoto, {\it Non-commutative spheres. II: rational rotations.} J. Operator Theory 27 (1992) 53--85.

\bibitem{casazza} P. G. Casazza and J. C. Tremain, {\it The Kadison-Singer problem in mathematics and engineering.} Proc. Nat. Acad. Sci. 103 no. 7 (2006) 2032--2039.

\bibitem{connes} A. Connes, {\it Classification of injective factors.} Annals of Mathematics 104 (1976) 73--115.

\bibitem{courtney1} K. E. Courtney, {\it C$^*$--algebras and their finite--dimensional representations.} PhD Dissertation, Department of Mathematics, University of Virginia, 2018.

\bibitem{courtney2} K. E. Courtney, {\it Kirchberg's QWEP conjecture:
between Connes' and Tsirelson's problems.} UK Operator Algebra Seminar, 2020.

\bibitem{courtney3} K. E. Courtney and T. Schulman, {\it Elements of $C^*$--algebras attaining their norm in a finite--dimensional representation.} Canadian J. Math. 71 no. 1 (2019) arXiv:1707.01949.

\bibitem{courtney4} K. E. Courtney and D. Sherman, {\it The universal $C^*$--algebra of a contraction.} J. Operator Theory 84 no. 1 (2020) 153--184. arXiv:1811.04043

\bibitem{courtney5} K. E. Courtney, {\it Universal $C^*$--algebras with the local lifting property.} Math. Scand. 127 (2021) 361--381. arXiv:2002.02365

\bibitem{davidson} K. R. Davidson, {\it $C^*$--Algebras by Examples.} American Math. Soc., Providence, Rhode Island, 1991.

\bibitem{dauns} J. Dauns, {\it The primitive ideal space of a $C^*$--algebra.} Canadian J. Math. Vol. XXVI, no. 1 (1974) 42--49.

\bibitem{elliot} G. A. Elliot and D. E. Evans, {\it The structure of irrational rotation $C^*$--algebras.} Annals of Mathematics 138 no. 3 (1993) 477--501.

\bibitem{gelfand} I. M. Gelfand and M. A. Naimark, {\it On the embedding of normed rings into the ring of operators on a Hilbert space.} Mat. Sbornik 12 no. 2 (1943) 197--217.

\bibitem{hall} B. Hall, {\it Lie Groups, Lie Algebras, and Representations, An Elementary Introduction.} Springer, Switzerland, 2003.

\bibitem{hatcher} A. Hatcher, {\it Vector Bundles and K--Theory.} https://openlibra.com/en/book/vector-bundles-and-k-theory

\bibitem{marcus} A. W. Marcus, D. A. Spielman, and N. Shrivastava, {\it Interlacing Famlies II: mixed characteristic polynomials and the Kadison--Singer problem.} Annals of Mathematics 182 no. 1 (2015) 327--350.

\bibitem{ji} Z. Ji, A. Natarajan, T. Vidick, J. Wright and H. Yuen, 
{\it MIP$^* = $ RE.} Communications of the ACM 64 no. 11 (2021) 131--138.
 arXiv:2001.04383.

\bibitem{kadison} R. V. Kadison and I. M. Singer, {\it Extensions of pure states.} American J. Math. 81 no. 2 (1959) 383--400.

\bibitem{lawton1} W. Lawton, A. Mouritzen, J. Wang, J. Gong {\it Spectral relationships between kicked Harper and on-resonance double kicked rotor operators.} J. Math. Phys. 50 no. 3 (2009) 032103. arxiv:0807.4276

\bibitem{lawton2} W. Lawton, {\it Minimal sequences and the Kadison--Singer problem}, Bull. Malaysian Math. Science Society 33 (2010) 169--176.

\bibitem{paulsen} V. Paulsen, {\it A dynamical system approach to the Kadison--Singer problem.} J. Functional Analysis 255 no. 1 (2008) 120--132.

\bibitem{neumann} J. von Neumann, {\it Allgemeine Eigenwerttheorie Hermitescher Funktionaloperatoren.} Math. Annalen 102 (1929) 49--131.

\bibitem{rammal} R. Rammal and J. Bellisard, {\it An algebraic semi-classical approach to Bloch electrons in a magnetic field.} J. Physics France 51 (1990) 1803--1830.

\bibitem{riesz} F. Riesz and B. Sz.-Nagy, 
{\it Functional Analysis.} Frederick Ungar Publishing Company, New York, 1955.

\bibitem{rudin} W. Rudin, {\it Real and Complex Analysis.} McGraw--Hill, Singapore, 1987.

\bibitem{shilov} G. Shilov, {\it Linear Algebra.} Dover, New York, 1977.

\bibitem{simon} B. Simon, {\it Almost periodic Schr$\ddot o$dinger operators: a review.}
Advances in Applied Math. 3 (1982) 463--490.

\bibitem{stacey} P. J. Stacey, {\it The automorphism groups of rational rotation algebras.}
J. Operator Theory 39 (1998) 395--400.

\bibitem{szigeti} J. Szigeti and L. Wyk, {\it A constructive elementary proof of the Skolem--Noether theorem for matrix algebras.} The American Mathematical Monthly 124 no. 16 (2017) 966--968. arXiv:1810.08368

\bibitem{waerden} B. van der Waerden, {\it Algebra I.} Springer, New York, 1991.

\bibitem{wang} H. Wang, D. Ho, W. Lawton, J. Wang, and J. Gong, {\it Kicked-Harper model versus on-resonance double-kicked rotor model: From spectral difference to topological equivalence.} Physical Review E 88 (2013) 1--15. 052920 arXiv:1306.6128

\bibitem{weaver} N. Weaver, {\it The Kadison--Singer problem in discrepancy theory.} 
Discrete Math. 278 no. 1--3 (2004) 227--239.

\bibitem{williams} D. Williams, {\it Crossed Products of $C^*$--Algebras.}
American Math. Society, 2007.

\bibitem{wintner} A. Wintner, {\it Zur Theorie der beschr$\ddot a$nkten Bilinearformen.} Math. Zeitschr. 30 (1929) 228--289.

\bibitem{yin} H. S. Yin, {\it A simple proof of the classification of rational rotation $C^*$--algebras.} Proceedings of the American Math. Society 98 no. 3 (1986) 469--470.

\end{thebibliography}
\end{document}